\documentclass{amsart}
\usepackage{amssymb,color}

\definecolor{c20}{rgb}{0.,0.7,0.}
\definecolor{c30}{rgb}{0.,0.,1.}
\definecolor{c40}{rgb}{1,0.1,0.7}
\definecolor{c50}{rgb}{1,0,0}

\def\crE#1{\textcolor{c50}{#1}}
\def\crEE#1{\textcolor{c30}{#1}}
\def\crEEE#1{\textcolor{c40}{#1}}

\def\crE#1{#1}
\def\crEE#1{#1}
\def\crEEE#1{#1}

\newtheorem{thm}{Theorem}
\newtheorem{cor}[thm]{Corollary}
\newtheorem{lem}[thm]{Lemma}

\theoremstyle{definition}

\newtheorem{exmp}{Example}
\newtheorem{remark}{Remark}

\newcommand{\BEX}{\begin{exmp}}
\newcommand{\EEX}{\end{exmp}}

\allowdisplaybreaks[4] 
  \def\Rset{\mathbb{R}}
  
\def\etal{et al.\ }

\newcommand{\nwc}{\newcommand}
\nwc{\COM}[1]{}

\newcommand{\nelem}[1]{{Lemma \ref{#1}}}

\newcommand{\netheo}[1]{{Theorem \ref{#1}}}

\newcommand{\kb}[1]{\boldsymbol{#1}}
\newcommand{\vk}[1]{\kb{#1}}

\newcommand{\ve}{\varepsilon}
\newcommand{\abs}[1]{\lvert #1 \rvert}
\newcommand{\Abs}[1]{ \Bigl \lvert #1 \Bigr \rvert}

\newcommand{\pk}[1]{\mbox{\rm$\vk{P}$} \{#1\} }

\newcommand{\pb}[1]{\mbox{\rm$\vk{P}$}\Bigl \{#1 \Bigr \}}

\def\R{\Rset}

\newcommand{\inr}{\in \R}

\newcommand{\ldot}{,\ldots,}

\newcommand{\limit}[1]{\lim_{#1 \to   \infty}}

\newcommand{\equaldis}{\stackrel{d}{=}}

\newcommand{\BQN}{\begin{eqnarray}}
\newcommand{\EQN}{\end{eqnarray}}
\newcommand{\BQNY}{\begin{eqnarray*}}
\newcommand{\EQNY}{\end{eqnarray*}}
\newcommand{\BS}{\begin{sat}}
\newcommand{\ES}{\end{sat}}
\newcommand{\BRM}{\begin{remark}}
\newcommand{\ERM}{\end{remark}}

\newcommand{\BL}{\begin{lem}}
\newcommand{\EL}{\end{lem}}
\newcommand{\BT}{\begin{thm}}
\newcommand{\ET}{\end{thm}}
\newcommand{\BK}{\begin{cor}}
\newcommand{\EK}{\end{cor}}

\newcommand{\QED}{\hfill $\Box$}
\newcommand{\IF}{\infty}

\def\fracl#1#2{\biggr( \frac{#1}{#2} \biggl) }

\newcommand{\prooftheo}[1]{ \textsc{Proof of Theorem} \ref{#1} }

\newcommand{\prooflem}[1]{\textsc{Proof of Lemma} \ref{#1}}

\def\Caro{ C_{a,\rho}}
\def\Aro{ \vk{\alpha}_{a,\rho}}

\begin{document}

\title[Tail asymptotics and estimation for elliptical distributions]{
Tail asymptotics and estimation for elliptical distributions}

\author{Enkelejd Hashorva}
\address{University of Bern,  Department of
Statistics\\
Sidlerstrasse 5, CH-3012 Bern, Switzerland}
\email{enkelejd.hashorva@stat.unibe.ch}

\subjclass[2000]{Primary 60F05; Secondary 60G70}

\keywords{Elliptical random vectors;
 exact asymptotics;  second order bounds; Gumbel max-domain of attraction; conditional  distribution; conditional
 quantile function;  joint survival probability; statistical estimation.}

\date{\today}
\begin{abstract}
Let $(X,Y)$ be a bivariate elliptical random vector with associated
random radius  in the Gumbel max-domain of attraction. In this paper
we obtain a second order asymptotic expansion of the joint survival
probability  $\pk{X > x, Y> y}$ for $x,y$ large.  Further, based on the asymptotic bounds we discuss
some aspects of the statistical modelling of joint survival
probabilities and the survival conditional excess probability.
\end{abstract}

\maketitle

\section{Introduction}
Let $(S_1,S_2)$ be a rotational invariant (spherical)  bivariate
random vector with associated random radius $R:=\sqrt{S_1^2+S_2^2}$.
The basic \crEE{distributional} properties of spherical random
vectors are obtained in Cambanis \etal (1981). So if $R>0$ almost
surely, then we have the stochastic representation
$$ (S_1,S_2) \equaldis (R U_1, R U_2),$$
where the \crEEE{bivariate} random vector  $(U_1,U_2)$ is
independent \crEEE{of the associated random radius} $R$ and
uniformly distributed on the unit circle of $\R^2$ ($\equaldis$
stands for equality of distribution functions). Linear combinations
of spherical random vectors define a larger class of random vectors,
namely that of elliptical random vectors. Canonical examples are 
the Gaussian and Kotz distributions (see Fang et al.\ (1990), Kotz et al.\ (2000), or Reiss and Thomas (2007)).\\
In this paper we consider a bivariate elliptical random vector defined in terms of $(S_1,S_2)$
and the pseudo-correlation coefficient \crE{$\rho\in (-1,1)$} via the stochastic
representation
\BQN \label{e}
 (X,Y)& \equaldis &(S_1, \rho S_1+ \sqrt{1- \rho^2} S_2).
 \EQN
Since for any $a,b$ two constants (see Lemma 6.2 of Berman
(1982))
\BQN\label{eq:ber} aS_1 + b S_2 \equaldis S_1 \sqrt{a^2+ b^2}
\EQN
we have $ X \equaldis Y \equaldis S_1.$ Referring to Cambanis et al.\ (1981) \crE{the random variable} $S_1$
is symmetric about $0$, and furthermore $ S_1^2 \equaldis R^2 W,$ with $W$ a Beta random variable with parameters $1/2, 1/2$
independent of $R$, implying that the distribution function of $X$ and $Y$
is completely known if the distribution function $F$ of $R$ is specified.

The basic distributional  properties of elliptical random vectors
are  well-known, see e.g., Kotz (1975), Cambanis et al.\ (1981),
Anderson and Fang (1990), Fang et.\ al (1990), Fang and Zhang
(1990), Berman (1992), Gupta and Varga (1993), Kano (1994), Szab{\l}owski (1998), or
Kotz et al.\ (2000)  among several others.

The tail asymptotic behaviour of each component, say of $X$,  can be
determined under some assumptions on the tail asymptotics of $R$.
The main work in this direction is done in Carnal (1970), Gale
(1980), Eddy and Gale (1981), Berman (1982, 1983, 1992) among
several others. For instance Berman (1992) shows the exact
asymptotic behaviour of $S_1$ if $R$ has distribution function in
the Gumbel max-domain of attraction. With motivation from Berman's
work, Hashorva (2007) obtains an exact asymptotic expansion of the
bivariate survival probability
$$ \pk{X> x, Y>ax}, \quad a\le 1$$
letting $x$ tend to $\IF$. See also Asimit and Jones (2007) for a partial result.

The main impetus for the present article comes from the recent deep
contribution  Abdous et al.\ (2007). We derive in this paper a refinement of
the asymptotic expansion of the joint survival probability obtained
in Hashorva (2007).  This is achieved assuming some second order asymptotic bounds
on the tail asymptotics of the distribution function $F$ as suggested in Abdous et al.\ (2007). \\

Our results are of a certain theoretical interest providing detailed
asymptotic expansions for a classical problem of probability theory
- tail asymptotics of random vectors. Further, based on our novel results, we suggest statistical estimators of the joint survival probability, conditional distribution and related quantile function.

We choose the following order for the rest of the paper: In Section
2 we present the main results.  Illustrating examples follow in Section 3. In Section 4 and Section 5 we discuss some
statistical aspects concerning the estimation of the joint survival
probability, conditional survival and conditional quantile functions.
Proofs and related asymptotics are relegated to Section 6.

\def\vax{v_a(x)}
\def\aroxy{\alpha_{\rho, x,y}}
\def\broxy{\beta_{\rho, x,y}}

\section{Asymptotic Bounds}

Let $(X,Y)$ be a bivariate elliptical random vector as in \eqref{e},
where the associated random radius $R$ has distribution function $F$
\crE{with upper endpoint $x_F\in (0,\IF]$ and} $F(0)=0$. Hashorva
(2007) derives an asymptotic expansion of the tail probability
$\pk{X> x, Y> ax},a\in (-\IF, 1]$ for $x\uparrow  x_F$ assuming that  $ F$ is in the
Gumbel max-domain of attraction, i.e., \BQN \label{eq:gumbel:w}
\lim_{u \uparrow x_F} \frac{1 - F(u+s/w(u))}{1- F(u)}&=& \exp(-s),
\quad \forall s\inr, \EQN with  $w$ some  positive scaling function.
Refer to Galambos (1987), Resnick (1987), Reiss (1989), Embrechts et
al.\ (1997),
 Falk et al.\ (2004), Kotz and Nadarajah (2005), or de Haan and Ferreira (2006) for
details on the Gumbel max-domain of attraction.\\
It is well-known (see e.g., Resnick (1987)) that \eqref{eq:gumbel:w}
is equivalent with the fact that the distribution function $F$ has
the following representation \BQN \label{eq:von2} 1- F(x)&=& d(x)
[1- F^*(x)], \quad x< x_F,
 \EQN
with $d(x)$ a positive function converging to $d\in (0,\IF)$ as $x\uparrow
x_F$ and \BQN\label{eq:von} F^*(x)&=& 1- \exp(- \int_{z_0}^x  w(s)
\,ds), \quad x< x_F, \EQN
 where $z_0$ is  a finite constant in the  left neighbourhood of the upper endpoint $x_F$. The distribution $F^*$ is referred to in the sequel as the von Mises distribution
related to $F$, whereas $w$ as the von Mises scaling function of
$F$.

Under the Gumbel max-domain of attraction assumption on $F$ (see
Hashorva (2007)) we have \BQN\label{eq:6}
 \pk{X> x, Y> ax} &=&  \crE{\frac{(1+o(1))c_1}{\crEEE{\sqrt{x w( \alpha_{\rho, x, ax} x)}}}}
 \pk{X> \alpha_{\rho, x,ax} x}, \quad x\to \IF,
 \EQN
\crEE{with $c_1$ a known constant}, provided that $a\in ( \rho, 1]$
and $x_F=\IF,$ where $\aroxy$ is defined by \BQN \label{eq:arho}
\aroxy&:=
\sqrt{1+(y/x- \rho)^2/(1- \rho^2)}. \EQN For $x,y$ positive such
that $y/x\ge a>\rho$ we have $ \aroxy > 1.$ Hence, since the scaling
function $w$ satisfies \BQN\label{eq:uv}
 \lim_{u \uparrow x_F } u w(u) &= &\IF, \quad
\text{ and }  \lim_{u \uparrow x_F } w(u) (x_F- u)  = \IF, \text{
if } x_F < \IF \EQN we conclude that the joint tail asymptotics in
\eqref{eq:6} is faster than
the componentwise tail asymptotics.\\
It turns out that for $y < \rho x$ or $y$ close enough to $\rho x$
the tail asymptotics of interest is up to a constant the same as
that of $\pk{X > x}$. \crE{For the bivariate Gaussian distribution
it is well-known that this fact is closely related to the so called
Savage
condition, see Dai  and Mukherjea (2001),  Hashorva and H\"usler (2003), or Hashorva (2005a) for details.}\\
For elliptical distributions the case where $y$ is close to $\rho x$
has been considered \crE{independently in Gale (1980),  Eddy and
Gale (1981)} and Berman (1982, 1983, 1992).

In this paper we are interested in a refinement of \eqref{eq:6}
which will be achieved under extra costs
related to a second order assumption on $F$.\\
Explicitly, as suggested in Abdous et al.\ (2007) we impose the
following assumption:

A1. {  \it Suppose that the distribution function $F$ satisfies
\eqref{eq:von2} with von Mises distribution function $F^*$, and
assume further that there exist positive measurable functions
$A_i,B_i,i=1,2$ such that for a scaling function $w$ for which
\eqref{eq:gumbel:w} holds we have \BQN\label{main:A} \Abs{ \frac{ 1-
F^*(u+ x/w(u))}{1- F^*(u)}- \exp(-x)} \le A_1(u)B_1(x) \EQN and
\BQN\label{main:A2} \Abs{ d(u+ x/w(u))- d(u)} & \le & A_2(u)B_2(x)
\EQN for all $u$ large and any $x\ge 0$. Furthermore, we assume that
$\limit{u} A_1(u)=\limit{u} A_2(u)=0,$ and $B_1,B_2$
 are locally bounded on any compact interval of $[0,\IF)$.
 }

Set in the following (whenever Assumption A1 is assumed) \BQN
\label{eq:Bs} A(x):= A_1(x)+ A_2(x), \quad B(x):=
B_1(x)+\exp(-x)B_2(x), \quad \forall x> 0. \EQN We note in passing
that $A_2(x)=B_2(x)=0,x>0$ is the original condition in the
aforementioned paper, where it is shown (see Lemma 7 therein) that
the class of distribution functions satisfying Assumption A1 is
quite large.

We consider in the sequel only distribution functions  $F$ with an
infinite upper endpoint. Further we assume that $\rho \in [0,1)$. We
state now the main result of this paper:

\BT \label{theo:1} Let $(S_1,S_2)$ be a bivariate spherical random
vector with associate random radius $R \sim F$, where the
distribution function $\crEEE{F}$ has an infinite upper endpoint and
$F(0)=0$. Let $(X,Y)$ be a bivariate elliptical random vector with
stochastic representation $(X,Y) \equaldis  (S_1, \rho S_1+ \sqrt{1-
\rho^2} S_2), \rho\in [0,1)$. Assume that \eqref{eq:gumbel:w} holds
with the scaling function $w$ and Assumption A1 is valid with
$A_i,B_i,i=1,2$ positive measurable functions
(with $A,B$ as in \eqref{eq:Bs}).\\
a) If $x,y$ are \crEEE{positive constants} such that \BQN
\label{eq:cond:mth:1} y \in (\rho x, x], \quad  \aroxy \ge c > 1,
\EQN and $\int_0^\IF B(s)\, ds < \IF$, then for all $x$ large \BQN
\label{eq:mthm:1} \pk{X> x, Y> y} &=& \frac{ \aroxy K_{\rho,x,y}}{2
\pi}\frac{1- F(\aroxy x)}{x w(\aroxy x)} \Biggl[1+ O\Biggl( A(\aroxy
x) +\frac{1}{x w(\aroxy x)}\Biggr)\Biggr], \EQN with $\aroxy$  as
defined in \eqref{eq:arho} and $K_{\rho,x,y}$ given by \BQN
\label{eq:kro} K_{\rho,x,y}&:= &\frac{x^2(1- \rho^2)^{3/2} }{(x-
\rho y)(y- \rho x)}\in (0,\IF). \EQN

b) Set $h(x):= xw(x),x>0$, and let $z_x,x\inr$ be given constants
such that $\abs{z_x}< K< \IF$ for all $x$ large.
 If further $\int_{ 0}^ \IF s^{-1/2}B(s) \, ds < \IF$, then for all $x$ large
and $y:= x \crE{[\rho  + z_x \sqrt{1- \rho^2}/ \sqrt{h(x)}+
\rho/h(x)]}$ we have \BQN\label{eq:m:2} \pk{X> x, Y> y    } &=&
\frac{1}{ \sqrt{2\pi}} \frac{ 1- F(x)}{ \sqrt{\crE{h(x)}}} [1-
\Phi(z_x)] \Biggl[1+ O\Biggl(
A(x)+\frac{1}{\crE{h(x)}}\Biggr)\Biggr], \EQN
\crEE{where $\Phi$ denotes the standard Gaussian distribution on $\R$.}\\
c) If $x,y$ are positive \crEEE{constants} with $y< a x, a\in
(0,\rho),\rho>0$ and $\int_{ 0}^ \IF s^{-1/2}B(s) \, ds < \IF$, then
\BQN \label{eq:m:3} \pk{X> x, Y> y  } &=& \frac{1}{ \sqrt{2\pi}}
\frac{ 1- F(x)}{ \sqrt{\crE{h(x)}}} \Biggl[1+ O\Biggl(
A(x)+\frac{1}{\sqrt{\crE{h(x)}}} \Bigl[  \frac{1}{\sqrt{\crE{h(x)}}}
+ \frac{1- F(\aroxy x)}{1- F(x)} \Bigr] \Biggr)\Biggr] \EQN is valid
for all large $x$. \ET

\begin{remark}
a) By the assumption on $F$ (recall \eqref{eq:uv}) we have
$ \limit{x}(x w(x))^{-1}= 0. $ If we assume that $\limit{u}A(u)= 0$ holds, then we have
$$  O\Biggl(   A(\aroxy x)
+\frac{1}{x w(\aroxy x)}\Biggr)= o(1), \quad x\to \IF.$$
Consequently by statement $a)$ in \netheo{theo:1} \BQNY \pk{X> x, Y>
y} &=& (1+o(1))\frac{\aroxy K_{\rho,x,y}}{2 \pi}\frac{1- F(\aroxy
x)}{x w(\aroxy x)}, \quad x\to \IF, \EQNY and if $y= a
x(1+o(1)),a\in (\rho, 1]$, then for all large $x$ we have
\BQNY 
\pk{X> x, Y> y} &=& (1+o(1))\frac{\alpha_{\rho,a}  K_{\rho,a}}{2 \pi}\frac{1-
F(\alpha_{\rho,a} x)}{x w(\alpha_{\rho,a}x)}, \quad x\to \IF, \EQNY holds with \BQN
\label{karo}
\alpha_{a, \rho}:= \sqrt{(1 -2a \rho+ a^2)/(1- \rho^2)}>1 ,\quad  K_{\rho,a}&:=&\frac{(1- \rho^2)^{3/2} }{(1- a\rho )(a- \rho )}\in (0,\IF).
 \EQN

b) Sufficient conditions for \eqref{main:A} to hold are derived in
Lemma 7  of Abdous et al.\ (2007). An instance is when $\limit{x}
w(tx)/w(x)=t^{\theta- 1}, \theta >0, \forall t>0$, i.e., the scaling
function $w$ defining $F^*$ is regularly varying at infinity with
index $\theta- 1$. Under the assumptions of Lemma 7 in the
aforementioned paper we may choose
\BQNY
A_1(u)= O(\abs{ (1/w(u))'}), \quad B_1(x)= (1+ x)^{- \kappa}, \quad u>0,x>0,\kappa >1.
\EQNY
c) We have $ \aroxy^2 = 1$ iff $y=\rho x$ and for any $y\in (0,x], x> 0$
$$ 1 \le \aroxy^2  \le \frac{2}{1+ \rho}, \quad \rho \in (-1,1).$$

d) Since \eqref{eq:gumbel:w} implies that  $1- F$ is rapidly varying
(see e.g., Resnick (1987) or de Haan and Ferreira (2006)) and
$\aroxy \ge C > 1$, then we have
$$ \limit{x} \frac{1- F(\aroxy x)}{1- F(x) }=0,$$
consequently in \eqref{eq:m:3} we have
\BQNY O\Biggl(
A(x)+\frac{1}{\sqrt{xw(x)}} \Bigl[  \frac{1}{\sqrt{xw(x)}} +
\frac{1- F(\aroxy x)}{1- F(x)} \Bigr] \Biggr)&=& O\Biggl(
A(x)+\frac{o(1)}{\sqrt{xw(x)}} \Biggr).
\EQNY
\end{remark}

We give next an alternative expansion of the tail probability under
consideration assuming  $y\in (\rho x, x]$.

\BT \label{theo:2} Let $(X,Y), \rho, F$ be as in \netheo{theo:1}.
Suppose that \eqref{eq:gumbel:w} holds with the scaling function $w$
and further Assumption A1 is satisfied with $A_i,B_i$ positive
measurable functions. Let $x,y$ be two positive constants such that
\eqref{eq:cond:mth:1} holds. If  $\int_0^\IF \max(1,
1/\sqrt{s})B(s)\, ds < \IF$ hods with $B$ as defined in
\eqref{eq:Bs},
 then for all $x$ large 
 we have
\BQN
\label{eq:mthm:2} 
\pk{X> x, Y> y}&=&  \frac{ \aroxy^{3/2} K_{\rho,x,y}}{\sqrt{2 \pi x
w(\aroxy x)}
} 
\pk{X > \aroxy x }\Biggl[1+ O\Biggl( A(\aroxy x) +\frac{1}{x
w(\aroxy x)}\Biggr)\Biggr],
 \EQN
 with   $\aroxy, K_{\rho,x,y}$ as defined in \eqref{eq:arho} and \eqref{eq:kro},
 respectively.
\ET \COM{ and \BQN \label{eq:Zcr} Z_{\zeta,\rho} &=& \frac{\zeta X+
Y}{\sqrt{\zeta^2+2\zeta\rho+ 1}} \equaldis X.}

Next, we consider the implications of our results for the tail asymptotics of the conditional survival probability $\pk{Y> y \lvert X> x}$.

\BT \label{theo:3} Under the assumptions and notation of \netheo{theo:2} we have for all $x,y$ large \BQN \label{eq:corr:2}
\pk{Y> y \lvert X> x} &=&   \frac{ \aroxy^{3/2} K_{\rho,x,y}}{\sqrt{2 \pi } } g_1(\aroxy, x)
= \frac{ \aroxy^{3/2} K_{\rho,x,y}}{\sqrt{2 \pi } } g_2(\aroxy, x),
 \EQN
 where
\BQN \label{eq:corr:g} g_1(\aroxy, x)&:= &\frac{\sqrt{w(x)/(\aroxy x)}}{
w(\aroxy x)}
 \frac{1- F(\aroxy x) }{ 1- F( x)} \Biggl[1+ O\Biggl( A(\aroxy x) +\frac{1}{x w(\aroxy x)}\Biggr)\Biggr],\\
 g_2(\aroxy, x)&:= &\frac{\pk{X> \aroxy x }}{ \sqrt{ x w(\aroxy x)} \pk{X>x}}\Biggl[1+ O\Biggl( A(\aroxy
x) +\frac{1}{x w(\aroxy x)}\Biggr)\Biggr].
\EQN
Furthermore, we have
\BQN \label{eq:limit} \limit{x} \pk{Y> y \lvert X> x} &=& 0,
\EQN
and for any $a\in (\rho, 1], z\ge 0$
\BQN \label{eq:corr:3} \pb{Y> ax+ z/w(\alpha_{a, \rho} x) \lvert X> x}&=&
\frac{ \alpha_{a, \rho}^{3/2} K_{a,\rho} }{\sqrt{2 \pi }} \exp(- \lambda_{a,\rho} z)  g_1(\alpha_{a, \rho}, x)\\
&=& \frac{ \alpha_{a, \rho}^{3/2} K_{a,\rho} }{\sqrt{2 \pi }} \exp(- \lambda_{a,\rho} z)  g_2(\alpha_{a, \rho}, x)
 \EQN
hold locally uniformly with respect to $z$,  with $\alpha_{a, \rho}, K_{a,\rho}$ as defined in \eqref{karo} and
\BQN\label{laro}
 \lambda_{\rho, a}&:=& \frac{ a - \rho}{ \sqrt{(1- \rho^2)(1 -2a \rho+ a^2)}}>0. \EQN \ET
Convergence to a positive constant in \eqref{eq:limit} is achieved
when $y$ depends on $x$ being close to $\rho x$ for $x$ large. See
Gale (1980), Berman (1982,1983,1992), Abdous et al.\ (2007),
Hashorva (2006,2007) for elliptical case, and Hashorva et al.\
(2007) for the more general case of Dirichlet distributions.

\begin{remark} a)   A slightly  more general setup is when Assumption A1 is reformulated considering
positive measurable functions $A_i, B_i, i\ge 1$ such that for all
$x>0$ and $u$ large \BQNY \Abs{ \frac{ 1- F^*(u+ x/w(u))}{1-
F^*(u)}- \exp(-x)} &\le & \sum_{i=1} ^\IF A_i(u) B_i(x)=:\Psi(u,x)
\EQNY is valid with $\limit{u} \sup_{i\ge 1} A_i(u)=0, i\ge 1,$
$\Psi(u,x)$ finite for all $x>0$ and $u$ large.

b) In the asymptotic results above the scaling function $w$ appears
prominently. One choice for the scaling function is the von Mises
scaling function, i.e., $w$ defines the von Mises distribution
function $F^*$ in \eqref{eq:von}. We can choose however another
scaling function $\overline{w}$ defined asymptotically by \BQN
\label{eq:w} \overline{w}(u)&:= & \frac{(1+o(1))[1-
F(u)]}{\int_{u}^{x_F} [1- F(s)]\, ds}, \quad u \uparrow  x_F. \EQN
Note that the bounding functions $A_i,B_i, i=1,2$ in the Assumption
A1 depend on which concrete scaling function we choose. In view of
Lemma 16 in Abdous et al.\ (2007) (see its assumptions and Lemma 7
therein) if $w_F$ is the von Mises scaling function defining $F^*$
in \eqref{eq:von} and $\overline{w}$ is another scaling function
defined by \eqref{eq:w}, then \eqref{main:A} holds with
$\overline{w}$ instead of $w_F$, and $\overline{A}_1$ instead of
$A_1$, where $\overline{A}_1 $ is defined by \BQN \label{eq:Ay}
 \overline{A}_1(u) &=& O \Bigl( \abs{(1/w_F(u))' } + \abs{w_F(u)- \overline{w}(u)}/ \overline{w}(u) \Bigr)
, \quad \forall u>0.
 \EQN
\end{remark}

If $F$ is a von Mises distribution function,  then we can make use
of
Lemma 7 in Abdous et al.\ (2007) (provided the assumptions of that lemma hold).\\
We discuss in the next lemma the case $F$ is  a mixture
distribution.

\BL \label{lem:klem} Let $F_i, i\ge 1$ be von Mises distribution
functions with the same scaling function $w$ and upper endpoint
infinity. Suppose that the Assumption A1 is satisfied for all
$F_i,i\ge 1$ with corresponding functions $w, A_{i1}, B_{i1}, i\ge
1$, and $A_{i2},B_{i2}, i\ge 1$ identical to 0. If $F$ is another
distribution function defined by
$$  F(x)= \sum_{i=1}^\IF a_i F_i(x), \quad \text{ with } a_i>0,i\ge 1\ :\sum_{i=1}^\IF a_i=1$$
for all $x$ large, then $F$ is in the Gumbel max-domain of
attraction. If $F$ is the distribution function of the random radius
$R$ in \netheo{theo:2} and \netheo{theo:3}, then both these theorems
hold with
$$ A(u)= \sum_{i=1} ^\IF A_{i1}(u), \quad u>0,$$
provided that $\sum_{i=1}^\IF a_i \int_0^\IF B_{i1}(s)\,ds < \IF$
and $A(u)$ is bounded for all large $u$ with $\limit{u} A(u)= 0$.
\EL

\section{Examples}
We present  next three illustrating examples.

{\bf Example 1.} [Kotz Type I]\\
Let $(X,Y)=R (U_1, \rho U_1+ \sqrt{1-\rho^2} U_2),\rho \in [0,1)$
with $R$ a positive random radius independent of the bivariate
random vector $(U_1,U_2)$ which is uniformly distributed on the unit
circle of $\R^2$. We call $(X,Y)$ a Kotz Type I elliptical random
vector if
 \BQNY
 1- F(x) &=&  C  x^{N}\exp(-c x^\delta), \quad c>0, C>0,\delta>0, N\inr
 \EQNY
 for any $x>0$. Set $w(u):= c \delta u ^{\delta-1},u>0$.  For any $x\inr$ we obtain
$$\limit{u} \frac{ \pk{R> u+ x/w(u)}}{\pk{R>u}} = \exp(-x).$$
Consequently, $F$ is in the Gumbel max-domain of attraction with the
scaling function $w$. Further, we have \BQNY \pk{X> u} &=& (1+o(1))
\frac{C}{\sqrt{2 r \delta \pi   }} u ^{N-\delta/2}\exp(-ru^\delta),
\quad u\to \IF. \EQNY The von Mises scaling function $w$ is given by
\BQNY w(u)
&=& c\delta  u^{\delta-1}   -  Nu^{\eta}, \quad \eta:=-1.
\EQNY Hence \eqref{eq:Ay} implies that the approximation under
consideration holds choosing $A(u):=A_1(u),u>0,$ where $ A_1(u):= O(
u^{-\delta}), u>0.$ Further we can take $B_1(x) = (1+ x)^{-\kappa},
\kappa>1$ and $B_2(x)=0, x> 0$. Hence the second order assumption in
our results above is satisfied. We may thus write for $x$ large and
$y \in (\rho x, x]$ such that $\aroxy > c>1$ \BQNY \pk{X> x, Y>  y}
\COM{ &=&(1+o(1))\frac{\Aro^{3/2} \Caro}{\sqrt{2 \pi}}
\fracl{1}{u_n  r \delta (\Aro u_n ) ^{\delta-1}  }^{1/2}\pk{X> \Aro u_n}\\
&=&(1+o(1))\frac{\Aro^{(3 - \delta)/2 +1} \Caro}{\sqrt{2 r \delta
\pi}} u_n  ^{-\delta/2}\pk{X> \Aro u_n}
\\
&=&(1+o(1))\frac{\Aro^{(3 - \delta)/2 +1/2} \Caro}{\sqrt{2 r \delta
\pi}}
u_n  ^{-\delta/2}\fracl{K}{2 r \delta \pi   }^{1/2} (\Aro u_n)  ^{N-\delta/2}\exp(-r(\Aro u_n)^\delta)\\
\\
&=&(1+o(1))\frac{\Aro^{(3 - \delta)/2 +1/2 + N-\delta/2}
\Caro}{\sqrt{2 r \delta  \pi}}
\fracl{K}{2 r \delta \pi   }^{1/2} u_n  ^{N}\exp(-r(\Aro u_n)^\delta)\\
\\
&=&(1+o(1))\frac{K \Aro^{(3 - 2\delta)/2 +N+1/2 } \Caro}{2 r \delta
\pi} u_n  ^{N-\delta}\exp(-r(\Aro u_n)^\delta).
\\
} &=&(1+o(1))\frac{C  \aroxy ^{2 - \delta +N } K_{\rho,x,y} }{2 r
\delta  \pi} x  ^{N-\delta}\exp(-r(\aroxy x)^\delta)
[1+ O( x^{-\delta})]\\
&=&(1+o(1))\frac{ \aroxy ^{2 - \delta/2 } K_{\rho,x,y} }{\sqrt{2 r
\delta  \pi}} x^{-\delta/2} \frac{C  \aroxy ^{N - \delta/2 }
K_{\rho,x,y} }{\sqrt{2 r \delta  \pi}}  x  ^{N-\delta}\exp(-r(\aroxy
x)^\delta)
[1+ O( x^{-\delta})]\\
&=&(1+o(1))\frac{ \aroxy ^{2 - \delta/2 } K_{\rho,x,y} }{\sqrt{2 r
\delta  \pi}} x^{-\delta/2}\pk{X>
\aroxy x}[1+ O( x^{-\delta})]\\
&=&(1+o(1))\frac{ \aroxy ^{2 - \delta/2 } K_{\rho,x,y} }{\sqrt{2 r
\delta  \pi}} \pk{X^*> \aroxy x}[1+ O( x^{-\delta})], \EQNY
with $(X^*,Y^*)$ another Kotz Type I random vector  with coefficients $C, N^*=N- \delta/2, r, \delta$.\\

{\bf Example 2.} [Tail Equivalent Distributions]\\
Let $G$ be a von Mises distribution function with scaling  function
$w$ such that \eqref{main:A} holds with functions $A_1,B_1$, and let
$F_{\gamma, \tau},\gamma, \tau>0$ be another distribution function
with infinite upper endpoint. Assume further that \BQN\label{gata}
 1- F_{\gamma, \tau}(x)= (1+ ax^{- \gamma}+ O(x^{- \tau \gamma})) [1- G(x)], \quad \forall x>0,
 \EQN
with $\gamma, \tau$ two positive constants. Clearly, $F$ and $G$ are
tail equivalent since
$$ \limit{x}  \frac{1- F_{\gamma, \tau}(x)}{1- G(x)}= 1.$$
Suppose further that $w(x):= c \delta x^{\delta- 1}, c>0, \delta >
0$, and set $d(x):=(1+ ax^{- \gamma}+ O(x^{- \tau -\gamma})), x> 0.$ We have for all large $u$ and $x>0$ \BQNY \lefteqn{\Abs{d(u+
x/w(u))- d(u)}}\\
&=&a u^{-\gamma} \Bigl[ (1+ x/(u w(u)))^{- \gamma} - 1\Bigr]- u^{-\tau -\gamma}\Bigl[ 1+ O((1+ x/(u w(u)))^{- \tau \gamma})\Bigr]\\
&=& A_2(u) B_2(x), \EQNY where $A_2(u):= O(u^{-[\gamma+ \min(\tau,
\delta)]}), u>0$  and $B_2(x)$ is such that $\int_0^\IF
B_2(x)\exp(-x)\, dx < \IF$. Hence our asymptotic results in the
above theorems hold for such $F$
with the function $A$ defined by $ A(u):= A_1(u)+ A_2(u), u>0$.\\

{\bf Example 3.} [Regularly Varying $w$]\\
\COM{Under the setup of the above example which leads to $w$ a power
function we consider $F$ a von Mises distribution function
so that 
$w$ is a regularly varying function. Explicitely,} Consider $F$ a
von Mises distribution function in the Gumbel max-domain of
attraction with the scaling function $w(x)= F'(x) /[1-F(x)],x>0$
defined by \BQNY
 w(x) &:=& \frac{c \delta x^{\delta -1}}{1+ t_1(x)} =
c \delta x^{\delta -1}(1+ o(1)), \quad c>0, \delta > 0, \quad
\forall x>0, \EQNY which implies (Abdous et al.\ (2007)) that all
$x$ large
$$  \pk{R> x} = \exp( - c x ^\delta(1+ t_2(x)))$$
holds, where  $t_i(x), i=1,2$ are two regularly varying  functions
with index $\eta \delta,\eta < 0$.

Now choosing $\overline{w}(x):=c \delta x^{\delta -1}$ instead of
$w(x)$ we conclude that the second order correction function $A_1$
satisfies
$$ A_1(u):= O( u^{-\delta}+ u^{\eta \delta} L_1(u)), \quad u> 0,$$
with $L_1(u)$ a positive slowly varying function, i.e., $\limit{u}
L_1(tu)/L_1(u)=1, t>0$. It can be easily checked that our main theorems above hold for this case with $A(u):=A_1(u), u>0$.\\

\COM{
{\bf Example 4.} [Infinite Sums]\\
Consider $F_i(x,\gamma_i, \tau_i,G_i), i\ge 1$ as specified in
\eqref{gata}, with $\gamma_1,\tau_i$ positive constants and $G_i,
i\ge 1$ von Mises distribution functions  as in the previous example
where the von scaling . Define a distribution function $F$ by
$$ F(x)= \sum_{i=1}^\IF a_i F_i(x, \gamma_i, \tau_i,G_i)(x), \quad x\ge 0, \quad a_i\ge 0, i\ge 1.$$
If we choose as scaling function $\overline{w}(x):= c \delta
x^{\delta- 1}, x>0$ then combining the two previous examples with
the result of \nelem{lem:klem}  we see that our asymptotic results
hold with this $w$ and correction function $A$ defined by
$$ A(u):=O( u^{- \min_{i \ge 1} [ \gamma_i+\min(\tau_i, \delta)]}+ u^{-\delta}+ u^{\eta \delta} L_1(u)),
\quad u> 0$$

}

\section{Estimation of Joint Survival Probability}
Consider the estimation of the joint survival probability $\pk{X>x, Y>y}$ for $x,y$ large,  with $(X,Y)$ a bivariate random vector satisfying the assumptions of \netheo{theo:2}. We discuss first the implications of \eqref{eq:mthm:1}, \eqref{eq:mthm:2}, and then outline the estimation motivated by \eqref{eq:m:2}.

\underline{Estimation based on \eqref{eq:mthm:1}, \eqref{eq:mthm:2}}\\
Consider the case $x=y$ is large. The constants $\alpha_{\rho, x,x},
K_{\rho, x,x}$ do not depend on $x$ and $y$ (given in \eqref{karo}
for $a=1,\rho\in [0,1)$). Both asymptotic expansions can be written
as
$$ \pk{X> x,Y> x}= q_1(x, \rho, w, F)= q_2(x, \rho, w, G_\zeta),$$
with $F$ and $G_\zeta$ the distribution function of $R$ and $Z_{\zeta, \rho}$, respectively.\\
If $(X_i,Y_i), 1 \le i\le n$ are independent  copies of $(X,Y)$, then a $\sqrt{n}-$consistent estimator of $\rho$ is available from the
literature. The more difficult part is the estimation of the tails
and the function $w$.\\
If we restrict ourselves to distribution functions $F$ in the Gumbel
max-domain of attraction with von Mises scaling function
$$w(x)= c \delta x^{\delta -1}, \quad c>0, \delta> 0, \quad x>0,$$
then $w$ can be estimated utilising the techniques in Abdous et al.\ (2007), where estimators for $c$ and $\delta$ are constructed using previous
results of  Girard (2004), Gardes and Girard (2006). See also the recent contribution Diebolt et al.\ (2007).\\
The estimation of the tail $1- F$ can be performed for instance if we restrict ourself to the case of Example 3 above.
Advanced extreme value statistics provide estimation of Gumbel tail and related quantiles under second order assumptions.
See the recent deep monographs de Haan and Ferreira (2006), Reiss and Thomas (2007). If we use \eqref{eq:mthm:2} then second order asymptotic
condition on the distribution of $X$ need to be imposed.\\
Note in passing that the scaling function $w$ appears in the assumption on $F$, which on the turn implies
that both  $X$ and $Y$ have distribution functions in the Gumbel max-domain of attraction
with the same scaling function $w$. Consequently we may estimate $w$
alternatively utilising only the observations $X_i,1 \le i \le n$, or more generally
we may estimate $w$ from the observations
$$Z_{\zeta,\rho,i}= \frac{\zeta X_i+ Y_i}{\sqrt{\zeta^2+2\zeta\rho+ 1}}\equaldis X, \quad \zeta \inr, \quad i=1 \ldot n.$$
Further, instead of estimating the tail $1- F$ we may estimate the tail $1-
G_\zeta$. If we use the random points $X_i, 1 \le i\le n$ to estimate $w$ the advantage is that the estimator of $\rho$ is not involved.
The disadvantage is that the second order correction is a consequence of an assumption on $R$ and not on $X$.

\underline{Estimation based on \eqref{eq:m:2}}\\
If $x,y$ are large positive constants then making use of the approximation in \eqref{eq:m:2} we may write (set $h(x):=\sqrt{xw(x)},x>0$)
$$ \pk{X>x,Y>y} \approx \frac{1}{\sqrt{2 \pi}}  \frac{1- F(x)}{\sqrt{h(x)}}
\Bigl[1- \Phi\Bigl( (y/x - \rho - \rho/h(x))\sqrt{h(x)/(1-
\rho^2)}\Bigr)     \Bigr],$$ hence an estimator of the probability
of interest can be constructing by the right hand side of the above
approximation. Again  we have the same estimation issues for the
tail $1- F$ and $w(x)$ where $x$ is large as above.

\section{Estimation of Conditional Survival and Quantile Function}
Let $(X,Y), (X_i,Y_i), 1 \le i\le n,$ be as in the previous section. Our asymptotic results
above can be employed for the estimation of the conditional excess survival function
$$\Psi(y,x):= \pk{Y> y \lvert X> x},$$
where $x,y$ are large. The estimation of the conditional distribution
$1- \Psi(y,x)$ is discussed in Abodus et al.\ (2007).\\
Clearly, one way to address this problem is to consider the estimation of the joint and marginal
survival probabilities $\pk{X> x, Y> y}, \pk{X>x}$, separately. As noted in the aforementioned paper
for large values of $x$ the empirical distribution function is useless since no observations might fall in the
relevant regions. \\
Our suggestion for the estimation of $\Psi(y,x)$ is motivated by the asymptotic relations
shown in \netheo{theo:3}. Under the assumptions of that theorem for all large $x$ and any $z\inr$ we have
\BQNY
\Psi(x+ z/w(\alpha_{1,\rho}x),x)&=& \frac{ \alpha_{1,\rho}^{3/2} K_{1,\rho} }{\sqrt{2 \pi }}
\exp(- \lambda_{1,\rho} z)  g_1(\alpha_{1,\rho}, x)\\
&=& \frac{ \alpha_{1,\rho}^{3/2} K_{1,\rho} }{\sqrt{2 \pi }} \exp(- \lambda_{1,\rho} z)  g_2(\alpha_{1,\rho}, x),
 \EQNY
with $\alpha_{1,\rho},K_{1,\rho}$ and $\lambda_{1,\rho}$ as in
\eqref{karo} and \eqref{laro}, respectively.\\
Let $\hat w_n, \hat g_{n1}, \hat g_{n2}, \hat \rho_n$ be estimators  of
the function $w, g_1,g_2$ and $\rho$, respectively, and set
$$\hat \alpha_n:=\alpha_{1,\hat \rho_n} , \quad \hat \lambda_n:= \lambda_{1,\hat \rho_n}, \quad n\ge 1.$$
Then  we may estimate $\Psi(y,x)$ by \BQN
\hat \Psi_{n,1}(y,x)&:=& \frac{ \hat \alpha_n^{3/2} K_{1,\hat \rho_n} }{\sqrt{2 \pi }} \exp(- \hat \lambda_n \hat w_n(\hat \alpha_n  x) (y- x))  \hat g_{n1}(\hat \alpha_n , x)\notag \\
&=:&\exp(- \hat \lambda_n \hat w_n(\hat \alpha_n  x) y)  \hat g^*_{n1}(\hat \alpha_n,x ),
\EQN
or alternatively
\BQN \hat
\Psi_{n,2}(y,x)&:=& \frac{ \hat \alpha_n ^{3/2} K_{1,\hat \rho_n}
}{\sqrt{2 \pi }} \exp(- \hat \lambda_n \hat w_n(\hat \alpha_n  x)
(y- x)) \hat g_{n2}(\hat \alpha_n , x)\notag \\
&=:& \exp(- \hat \lambda_n \hat w_n(\hat \alpha_n
x) y) \hat g^*_{n2}(\hat \alpha_n , x).
\EQN
Let $\Phi(q,x), q\in (0,1),x>0$ be the quantile function defined as the inverse
function of $1- \Psi(y,x)$ for $x>0$ fix. Inverting the above
expressions we have also two estimators of the conditional quantile function
\BQN \hat y_{n,1}(q,x)&:=& \frac{\ln ( \hat
g^*_{n1}(\hat \alpha_n, x ))- \ln (1- q)}{ \hat \lambda_n
\hat w_n(\hat \alpha_n  x)},
\EQN and
\BQN \hat y_{n,1}(q,x)&:=& \frac{\ln ( \hat g^*_{n2}(\hat \alpha_n , x ))- \ln
(1- q)}{\hat \lambda_n \hat w_n(\hat \alpha_n
x)}.
\EQN
The difficulty in constructing these estimates lies in the fact that we have to estimate both the
tail  $1- F$ (implicit in the estimation of the functions $g_1,g_2$), and the scaling function  $w$.
In the setup of Example 3 above there is a simple relation about these
functions, and the second order correction is easy to handle.\\

\section{Proofs and Related Asymptotics}
In the following lemma we derive a formula for the distribution
function of a bivariate elliptical random vector. Define next for
$a\ge  1$ and $x,y$ positive constants \BQN\label{iax} I(a,x):=
\int_{a}^\IF [1- F(xs)]\frac{1}{s\sqrt{s^2- 1}}\,ds, \EQN and
\BQN\label{eq:arobro}
 \aroxy:= \sqrt{1+ ((y/x)- \rho)^2/(1- \rho^2)}\ge 1,\quad
\crE{\broxy}:= \aroxy x/y, \quad x,y\inr, \rho\in(-1,1).
\EQN
\BL
Let $(S_1,S_2)$ be a bivariate spherical random vector with associate random radius $R$ which has distribution function $F$.
If $\rho \in \crEEE{[0,1)}$  and $F(0)=0$ then we have:
a) If \crEEE{$x>0$ and $y\in (\rho x, x]$}
\BQN \label{eq:lem:A:A}
 \pk{\crE{S_1 > x, \rho S_1+ \sqrt{1- \rho^2} S_2}> y}&=& \frac{1}{2 \pi} [I( \aroxy,x)+ I(\crE{\broxy},y)].
\EQN
b) If $y/x < \rho$ and $x>0, y\ge 0$
\BQN\label{eq:lem:A:B}
 \pk{\crE{S_1 > x, \rho S_1+ \sqrt{1- \rho^2} S_2}> y}&=& \crE{\frac{1}{2 \pi}}[\crE{2}I(1,x) - I( \aroxy,x)+ I(\crE{\broxy},y)].
\EQN
\EL
\begin{proof} We have the stochastic representation (see Cambanis et al.\ (1981), Berman (1992))
\BQNY
(S_1, \rho S_1+ \sqrt{1- \rho^2} S_2) &=& (S_2, \rho S_2+ \sqrt{1- \rho^2} S_1)\\
 &\equaldis& (S_1 \cos(\Theta), S_2 \cos(\Theta - \psi)),
\EQNY
with $\Theta$ uniformly distributed in $(-\pi, \pi)$ independent of $R$ and $\psi:= \arccos(\rho)$.
For $x>0,\crEEE{y}\ge 0$  two constants we may thus write (see Lemma 3.3 in Hashorva (2005b))
\COM{
\BQNY
\pk{S_1> x, \rho S_1+ \sqrt{1- \rho^2} S_2> y}&=&
\pk{R \cos (\Theta) > x, R \crE{\cos (\Theta- \psi)} > y}.
\EQNY
}
\BQNY
\lefteqn{
2 \pi \pk{R \cos (\Theta) > x, R \crE{\cos (\Theta- \psi)}> y}}\\
&=&
\int_{\arctan((y/x- \rho)/\sqrt{1- \rho^2})}^{\pi/2} \pk{ R> x/ \cos(\theta)} \, d \theta
\crE{+ \int_{\arccos(\rho)- \pi/2}^{\arctan((y/x- \rho)/\sqrt{1- \rho^2})} \pk{ R> \crEEE{y}/ \cos(\theta- \psi))}
\, d \theta}\\
&=&\int_{\arctan((y/x- \rho)/\sqrt{1- \rho^2})}^{\pi/2} \pk{ R> x/ \cos(\theta)} \, d \theta+ \int_{\arctan((y/x- \rho)/\sqrt{1- \rho^2})}^{ \pi/2 -\psi}
\pk{ R> \crEEE{y}/ \cos(\theta+ \psi) }\, d \theta.
\EQNY
As in Abdous et al.\ (2007) we obtain for  $y/x \ge \rho$
\BQNY
2 \pi \pk{R \cos (\Theta) > x, R \crE{\cos (\Theta- \psi) }> y}&=&
I( \aroxy,x)+ I(\crE{\broxy},y),
\EQNY
and if $y/x < \rho$ with $x,y$ positive
\BQNY
2 \pi \pk{R \cos (\Theta) > x, R \crE{\cos (\Theta- \psi)}> y}&=&
\crE{2}I(1,x) - I( \aroxy,x)+ I(\crE{\broxy},y),
\EQNY
with $\aroxy, \broxy$ as \crEEE{defined} in \eqref{eq:arobro}, hence the proof is complete. \end{proof}

Note in passing that
$$(S_1 \cos(\Theta), S_2 \cos(\Theta - \arccos(\rho)))\equaldis  (S_1 \cos(\Theta), S_2 \sin(\Theta + \arcsin(\rho))),$$
which leads to the alternative formula derived in Abdous et al.\
(2007) and Kl\"uppelberg et al.\ (2007) for the tails of elliptical
distributions. Remark further that some alternative formulae for the
distribution of bivariate elliptical random are presented in Lemma
3.3 in Hashorva (2005b).

\crE{In the next lemma we consider a real function $a(x)> 1, \forall x\inr$. We write for notational simplicity $a$ instead of $a(x)$.}

\BL \label{lem:1} Let $F$ satisfy \eqref{eq:gumbel:w} with the scaling function $w$ and $x_F=\IF, F(0)=0$.
Assume further that the Assumption A1 holds.\\
i) If  $\int_0^\IF B(s)\, ds < \IF$, then for any function $a:=a(x)> 1, x\inr$ and $x$ large we have
\BQN\label{lem:1:a}
\int_{a}^\IF [1- F(xs)] \frac{1}{ s\sqrt{s^2- 1}} \,ds &=&  \frac{1}{ a\sqrt{a^2- 1}}
\frac{1- F(ax)}{x w(ax)} \Biggl[1+ O\Biggl(  A(ax)+\frac{1}{xw(ax)}\Biggr)\Biggr].
\EQN
ii) If $z_x,x\inr$ is such that for all $x$ large $0 \le z_x< K< \IF $ and further
$ \int_{z_x}^\IF B(s)/\sqrt{s} \, ds < \IF$, then for all large $x$ we have
\BQN\label{lem:1:b}
\int_{1+ z_x/(xw(x))}^\IF [1- F(xs)] \frac{1}{ s\sqrt{s^2- 1}} \,ds &=&
\frac{1- F(x)}{\sqrt{xw(x)}} \sqrt{2\pi}  [1- \Phi(\sqrt{2z_x})] \Biggl[1 +   O\Biggl(  A(x)+\frac{1}{xw(x)}\Biggr)\Biggr].
\EQN
\EL
\begin{proof}
i) Let $x$ be a given positive constant. \crEE{Set
$$a:=a(s), \quad v_a(s):=s w(as), \quad s\ge 0$$
and define}
$I(a,x), \aroxy, \broxy$ as in \eqref{iax} and \eqref{eq:arobro}, respectively. Transforming the variables we have
\BQNY
I(a,x)
&=&\frac{1- F(ax)}{\vax} \int_{0}^\IF \frac{1- F(ax+ s/w(ax)) }{1- F(ax)} \frac{1}{(a+ s/\vax )\sqrt{(a+ s/\vax )^2- 1}}\,ds.
\EQNY
Further \eqref{main:A} implies for any $s\ge 0$ and all $x$ large
$$ \Abs{\frac{1- F(ax+ s/w(ax)) }{1- F(ax)}  - \exp(- s)} \le A(ax) B(s).$$
\COM{
Further for any $s\ge 0$
$$ \frac{1}{a}-\frac{1}{a+ s/\vax } = \frac{1}{\vax}\frac{s}{a(a+ s/\vax )}\le
\frac{1}{\vax}\frac{s}{a^2}$$
\BQNY
 \frac{1}{\sqrt{a^2-1}} - \frac{1}{ \sqrt{(a+ s/\vax )^2- 1}}
 &\le &  \frac{ 2 a  s/\vax+ (s/\vax)^2}{ 2(a^2-1)}.
\EQNY
The last inequality follows easily since for any $b>0, s\ge 0$
$$ \frac{1}{\sqrt{b}}- \frac{1}{\sqrt{b+ s}}= \frac{s}{\sqrt{b+ s}\sqrt{b}(\sqrt{b+ s}+ \sqrt{b})}
\le  \frac{s}{2 b}, \quad \forall s \ge 0.$$

For $f,g,f_1,g_1$ four real valued functions we may write 
\BQN\label{eq:simple} \abs{f(s)g(s)- f_1(s)g_1(s)} \le \abs{f(s)-
f_1(s)}\abs{g(s}+ \abs{g(s)- g_1(s)}\abs{f_1(s)}, \quad s\inr. \EQN
Applying the above we obtain for any $s\ge 0$ \BQNY
\Abs{\frac{1}{(a+ s/\vax )\sqrt{(a+ s/\vax )^2- 1}}   -  \frac{1}
{a\sqrt{a^2- 1}}}
 &\le &
\frac{1}{\vax}\Bigl[\frac{1}{a^2\sqrt{a^2-1}} + \frac{ 1  }{
a^2-1}\Bigr]  s, \EQNY }
Hence utilising further the Assumption A1 we may write for all $x$ large 
\BQNY\lefteqn{ \Abs{\frac{\vax} {1- F(ax)}\int_{a}^\IF [1- F(xs)]\frac{1}{s\sqrt{s^2- 1}}\,ds - \int_0^\IF \frac{\exp(-s)}
{a\sqrt{a^2- 1}} \, ds }}\\
&\le & \int_{0}^\IF
\Abs{\frac{1- F(ax+ s/w(ax)) }{1- F(ax)}  -\exp(-s)} \frac{1}{(a+ s/\vax )\sqrt{(a+ s/\vax )^2- 1}}\, ds\\
&&+ \int_{0}^\IF \exp(-s) \Abs{\frac{1}{(a+ s/\vax )\sqrt{(a+ s/\vax )^2- 1}}   -  \frac{1}
{a\sqrt{a^2- 1}}}\, ds\\
&\le &  \frac{1}{a\sqrt{a^2- 1}} A(ax) \int_{0}^\IF  B(s)\, ds+O( \frac{1}{\vax}) \\
&= &  O( A(ax)+ \frac{1}{\vax}),
\EQNY
thus the first claim follows.\\

\def\vax{\crE{h(x)}}

ii) Next we consider the second case $a(x):=1+ z_x/(xw(x)), x>0$. Set $\vax:=xw(x),x>0$.
Transforming the variables we obtain for all $x,y$ positive
\BQNY
 I(1+ z_x/\vax,x)&=&
 \frac{1- F(x)}{\sqrt{\vax}} \int_{z_x}^\IF \frac{1- F(x+ s/w(x)) }{1- F(x)} \frac{1}{\sqrt{\vax}(1+ s/\vax )\sqrt{(1+ s/\vax )^2- 1}}\,ds.
\EQNY
\COM{
For any $s\ge 0$ \eqref{eq:simple} implies
\crEEE{
\BQNY
\Abs{\frac{1}{\sqrt{2s}} - \frac{1}{ (1+ s/\vax) \sqrt{\vax(1+ s/\vax)^2- \vax}} }
&\le & 5/4 \sqrt{s}/(\vax \sqrt{2}).
\EQNY
}
}
%
\COM{Note that for $s$ and $x$ positive we have
$$1 - \frac{1}{1+ s/\vax }= \frac{1}{ \vax} \frac{s}{1+ s/\vax }\le \frac{1}{\vax} s,$$
$$\frac{1}{ \sqrt{(1+ s/\vax )^2- 1}} =
\frac{1}{ \sqrt{1 + 2s/\vax + (s / \vax)^2} - 1}
= \frac{1}{ \sqrt{s/\vax + (s / \vax)^2/2} }
= \frac{\sqrt{\vax}}{ \sqrt{2s+ s^2 / \vax} }, $$
and
$$ \frac{1}{\sqrt{2s}}- \frac{1}{ \sqrt{2s+ s^2 / \vax} }
\le \frac{  s^2/ \vax}{2 2 s} = \frac{1}{ 4 \vax} s.$$
}
\COM{
\crEEE{Note that for any  $s\ge 0$
\BQNY
\frac{1}{\sqrt{2s}} - \frac{1}{ (1+ s/\vax) \sqrt{\vax(1+ s/\vax)^2- \vax}}  &= &
\frac{1}{\sqrt{2s}}- \frac{1}{ (1+ s/\vax) \sqrt{2s+ s^2/\vax}}  \\
&= &\frac{1}{\sqrt{2s}}\Abs{1-  \frac{1}{1+ s/\vax}}\\
&&+ \frac{1}{1+ s/\vax}\Abs{ \frac{1}{\sqrt{2s}} - \frac{1}{ \sqrt{2s+ s^2/\vax}}}\\
&\le &\frac{1}{\sqrt{2s}}s/\vax +\frac{\sqrt{2s+ s^2/\vax}- \sqrt{2s} }{\sqrt{2s} \sqrt{2s+ s^2/\vax}}
\\
&=&\frac{1}{\sqrt{2s}}s/\vax + \frac{s^2/\vax}{(\sqrt{2s+ s^2/t_n}+ \sqrt{2s} ) \sqrt{2s} \sqrt{2s+ s^2/\vax}}\\
&=& \frac{1}{\sqrt{2s}}s/\vax + \frac{s^2/\vax}{2 \sqrt{2s}  \sqrt{2s} \sqrt{2s}}\\
&=& \frac{1}{\sqrt{2s}} (s/\vax + s/(4\vax)) = 5/4 \sqrt{s}/(\vax
\sqrt{2}). \EQNY } } Hence since $\int_{z_x}^\IF s^{-1/2}B(s) \,ds <
\IF$ for all $x$ large we have (recall Assumption A1)
\BQNY \lefteqn{
\Abs{ \int_{z_x}^\IF \frac{1- F(x+ s/w(x)) }{1- F(x)} \frac{
1}{\sqrt{\vax}(1+ s/\vax )\sqrt{(1+ s/\vax )^2- 1}}\,ds
-  \int_{z_x}^\IF \exp(-s) \frac{ 1 }{\sqrt{2s}}\,ds }}  \\
&\le &
\int_{z_x}^\IF \Abs{ \frac{1- F(x+ s/w(x)) }{1- F(x)} - \exp(-s)}\frac{1}{\sqrt{\vax}(1+ s/\vax )\sqrt{(1+ s/\vax )^2- 1}} \,ds
\\
&& +  \int_{z_x}^\IF \exp(-s)\Abs{ \frac{1}{\sqrt{\vax}(1+ s/\vax )\sqrt{(1+ s/\vax )^2- 1}}-
\frac{ 1 }{\sqrt{2s}}} \,ds   \\
&\le &
A(x) \int_{z_x}^\IF B(s)\frac{1}{\sqrt{2s} } \,ds
+  O(\frac{1}{\vax} )\int_{0}^\IF  \exp(-s) \sqrt{s}\,ds\\
&=&O(A(x)+ \frac{1}{\vax}).
\EQNY
Consequently \crE{since $z_x$ is bounded for all $x$ large }
\BQNY
I(1+ z_x/\vax,x) & =&  \frac{1- F(x)}{\sqrt{\vax}} \int_{z_x}^\IF
\exp(-s) \frac{1}{\sqrt{2s}}\,ds \Biggl[ 1+   O(A(x)+ \frac{1}{\vax})\Biggr]\\
&=&\frac{1- F(x)}{\sqrt{\vax}} \sqrt{2 \pi } [1- \Phi(\sqrt{2z_x})] \Biggl[ 1+   O(A(x)+ \frac{1}{\vax})\Biggr]
\EQNY
is valid with \crE{$\Phi$ the standard Gaussian distribution function on $\R$.}
 Thus the proof is complete. \end{proof}


\def\pxy{p_{xy}}
\prooftheo{theo:1} \crE{Define $I(a,x)$ and $\aroxy, \broxy$ as in
\eqref{iax} and \eqref{eq:arobro}, respectively.} Assumption A1
implies  that for any $x>0, \ve >0$ and  $u$ large we have \BQNY
\lefteqn{\Abs{ \frac{ 1- F(u+ x/w(u))}{1- F(u)}- \exp(-x)}}\\
&\le & \frac{d(u+ x/w(u))}{d(u)}  \Abs{ \frac{ 1- F^*(u+ x/w(u))}{1-
F^*(u)}- \exp(-x)  }
+ \exp(-x) \Abs{ \frac{ d(u+ x/w(u))} { d(u)} -1 } \\
  &\le & (1+ \ve) A_1(u)B_1(x) +\exp(-x)A_2(u) B_2(x),
\EQNY with $F^*$ the von Mises distribution function given in
\eqref{eq:von}. We assume for simplicity in the following that the
function $d(\cdot)$ is a constant for all $x$ large, impying
$A_2(u)=0$ for all $u$ large.

a) Let $x,y$ be two positive constants such that $y> \rho x$ and
$\aroxy \ge c>1$. In order to complete the proof we need a formula
for the survival probability \crEE{$\pk{X> x, Y> y}$}. In view of
\eqref{eq:lem:A:A} we have \COM{ \BQNY \pk{X> x, Y > y} =
\frac{1}{2\pi} \int_{ \arctan ( y/x- \rho)/(1 - \rho)}^{\pi/2} [1-
F()]\, d \theta + \int_{- \arcsin(\rho)}^{} \EQNY Since } for  $y>
\rho x$ and $x,y$ positive \BQNY 2 \crE{\pi} \crEE{\pk{X> x, Y> y}}
&=& I(\aroxy,x)+ I(\broxy,y). \EQNY Since further $\aroxy \ge c>
1$, applying \nelem{lem:1} we obtain \BQNY
\lefteqn{2 \pi \crEE{\pk{X> x, Y> y}}}\\
&=& \frac{1}{ \aroxy \sqrt{\aroxy ^2- 1}}
\frac{1- F(\aroxy x)}{x w(\aroxy x)} \Biggl[1+ O\Biggl( A(\aroxy x)+\frac{1}{xw(\aroxy x)}\Biggr)\Biggr]\\
&&+\frac{1}{ \aroxy (x/y)\sqrt{\aroxy ^2(x/y)^2- 1}} \frac{1-
F(\aroxy x)}{y w(\aroxy x)}
 \Biggl[1+ O\Biggl(   A(\aroxy x)+\frac{1}{xw(\aroxy x)}\Biggr)\Biggr]\\
&=&\frac{1- F(\aroxy x)}{x w(\aroxy x)} \Biggl[\frac{1}{ \aroxy
\sqrt{\aroxy ^2- 1}}+ \frac{1}{ \aroxy \sqrt{\aroxy ^2(x/y)^2- 1}}
+ O\Biggl(   A(\aroxy x)+ \frac{1}{xw(\aroxy x)}\Biggr)\Biggr]\\
&=&\crE{\frac{1- F(\aroxy x)}{x w(\aroxy x)} \Biggl[ \frac{x^2 (1-
\rho^2)^{3/2}\aroxy}{(y- \rho x)(x- \rho y)}
+ O\Biggl(   A(\aroxy x)+ \frac{1}{xw(\aroxy x)}\Biggr)\Biggr]}.
\EQNY

b) Let $z_x,x\inr$ be constants bounded for all $x$. We may write
for all $x,y$
$$\crEE{\pk{Y> y \lvert X> x}}= \pk{X> x}\frac{ \pk{X> x, Y> y}}{\pk{X> x}} =:  \pk{X> x}\chi(x,y).$$
For any $x$ positive 
$$ \pk{X> x} = 
\crEE{\frac{1}{\pi}}I(1,x),$$ hence \nelem{lem:1} implies for all
$x$ large enough
\BQNY
\pk{X> x} &=& \frac{1}{ \sqrt{2\pi}} \frac{ 1- F(x)}{ \sqrt{h(x)}}
\Biggl[1+ O\Biggl(   A( x)+\frac{1}{h(x)}\Biggr)\Biggr], \EQNY
\crE{with $h(x):=xw(x),x>0$}. In view of Theorem 3 in Abdous et al.\
(2007) we have for all large $x$ \BQNY \chi(x,y) &=& [1-
\Phi(z_x)]\Biggl[1+ O\Biggl( A(x)+\frac{1}{h(x)}\Biggr)\Biggr].
\EQNY Hence for all large $x$ \BQNY \crEE{\pk{X> x, Y> y }} &=
&\frac{1}{ \sqrt{2\pi}} \frac{ 1- F(x)}{ \sqrt{h(x)}} [1- \Phi(z_x)]
\Biggl[1+ O\Biggl( A(x)+\frac{1}{h(x)}\Biggr)\Biggr]^2, \EQNY
\crEE{thus the result follows. }\\
c) Now we consider the last case. By \eqref{eq:lem:A:B} for all
$x,y$ large $y \le ax < \rho x$ (implying $y \in (0, x]$)
 \BQNY \crE{2 \pi} \crEE{\pk{X> x, Y> y}}&= &2
I(1,x)- I(\aroxy ,x)+ I(\broxy ,y), \EQNY where $\aroxy\ge   \sqrt{
1+ (a- \rho)^2/(1- \rho)^2} >  1 $ and $\broxy:=\aroxy x/y>1$. By
the above results we have for all $y < ax$ and $x,y$ large enough
\COM{\BQNY I(\aroxy ,x)+ I(\broxy ,y)&=&\frac{1- F(\aroxy x)}{2 \pi
x w(\aroxy x)} \Biggl[1+ O\Biggl(
A(x)+\frac{1}{h(x)}\Biggr)\Biggr]^2. \EQNY Using further
\eqref{eq:uniG} we obtain for $x$ large } \BQNY \crEE{\pk{Y> y
\lvert X> x}} &=& \frac{1}{ \sqrt{2\pi}} \frac{ 1- F(x)}{
\sqrt{h(x)}}
\Biggl[1+ O\Biggl(   A(x)+\frac{1}{h(x)} + \frac{1}{\sqrt{h(x)}}
 \frac{1- F(\aroxy x)}{1- F(x)} \Biggr)\Biggr],
\EQNY thus the proof is complete. \QED

\prooftheo{theo:2} From the proof of \netheo{theo:1} we obtain for
all large $x$ \BQNY
 \pk{X>x}&=&\frac{1- F(x)}{\sqrt{2 \pi x w(x)}} \Biggl[1+ O\Biggl( A(x) +\frac{1}{x w(x)}\Biggr)\Biggr].
 \EQNY
Since further $\aroxy\ge 1$ and it is bounded for all $x$ large and $y$ positive such that (12) holds,
we may write for all $x$ large 
\BQN \label{eq:xxxxA}
 \pk{X> \aroxy x}&=&\frac{1- F(\aroxy x)}{\sqrt{2 \pi \aroxy x w(\aroxy x )}}
  \Biggl[1+ O\Biggl( A(\aroxy x) +\frac{1}{x w(\aroxy x)}\Biggr)\Biggr].
 \EQN
Applying \netheo{theo:1} we have
\BQNY
\lefteqn{\pk{X> x, Y> y}}\\
&=& \frac{ \aroxy^{3/2} K_{\rho,x,y}}{\sqrt{2 \pi x w(\aroxy x)}
\pi}\frac{1- F(\aroxy x)}{\pk{X> \aroxy x} \sqrt{ 2 \pi \aroxy  x w(\aroxy x)}} \pk{X> \aroxy x}\\
&& \times \Biggl[1+ O\Biggl( A(\aroxy
x) +\frac{1}{x w(\aroxy x)}\Biggr)\Biggr]\\
&=&\frac{ \aroxy^{3/2} K_{\rho,x,y}}{\sqrt{2 \pi x w(\aroxy x)}}\pk{X> \aroxy x}\Biggl[1+ O\Biggl( A(\aroxy
x) +\frac{1}{x w(\aroxy x)}\Biggr)\Biggr]^2\\
&=&\frac{ \aroxy^{3/2} K_{\rho,x,y}}{\sqrt{2 \pi x w(\aroxy x)}
}\pk{X> \aroxy x}\Biggl[1+ O\Biggl( A(\aroxy
x) +\frac{1}{x w(\aroxy x)}\Biggr)\Biggr].
\EQNY
Thus the result follows. \QED
\COM{In view of \eqref{eq:ber} we have for any $\zeta\inr$
\BQNY
 \zeta X + Y 
 &=& ( \zeta + \rho)S_1+ \sqrt{1- \rho^2} S_2\\
&\equaldis &\sqrt{( \zeta +\rho^2)+ 1- \rho ^2} S_1\\
&\equaldis &\sqrt{\zeta^2+2\zeta\rho+ 1} X
\EQNY
which yields
$$ \pk{X> z}= \pk{  \zeta X + Y> \sqrt{( \zeta +\rho^2)+ 1- \rho ^2}  z}, \quad \forall z>0,$$
thus the proof is complete. \QED
}
\COM{
\BQNY
\lefteqn{\pk{X> x, Y> y}}\\
 &=&\frac{ \aroxy^{3/2} K_{\rho,x,y}}{\sqrt{2 \pi x w(\aroxy x)}
}\pb{\frac{\zeta X+ Y}{\sqrt{(c+\rho^2)+ 1- \rho ^2}}> \aroxy x }\Biggl[1+ O\Biggl( A(\aroxy
x) +\frac{1}{x w(\aroxy x)}\Biggr)\Biggr].
\EQNY
}

\prooftheo{theo:3} The proof of the first claim follows easily since
$$\pk{Y> x \lvert X> x}= \pk{X>x,Y> y}/\pk{X>x}, \quad \forall x>0$$
utilising further the results of \netheo{theo:1} and \eqref{eq:xxxxA}.
The fact that the distribution function $F$ is in the Gumbel max-domain of attraction with scaling function $w$
implies that also the random variable $X$  has distribution function in the Gumbel max-domain of attraction (Berman (1992), Hashorva (2005b)) with the same scaling function $w$.
If  $y$ are positive constants such that $\aroxy> c> 1$ holds, then for all large $x$
$$ \limit{x} \frac{\pk{X> \aroxy x }}{\pk{X>x}} =0,$$
hence the result follows.\QED

\prooflem{lem:klem} Since $F_i,i\ge 1$ are in the Gumbel max-domain
of attraction with the same scaling function $w$ if follows easily
that the distribution function $F$ is in the Gumbel max-domain of
attraction with the same scaling function $w$. It can be easily
checked that both \netheo{theo:2} and \netheo{theo:3} hold with $
A(u):= \sum_{i=1} ^\IF A_{i1}(u), u>0$, hence the result follows.
\QED

\end{document}